\documentclass[11pt,twoside,leqno]{aomamlt2e} 
  \pageno{1353}
\received{November 24, 2003}

     \newcommand{\BA}{{\mathbb {A}}} 
     \newcommand{\BC}{{\mathbb {C}}}

     \newcommand{\BQ}{{\mathbb {Q}}} \newcommand{\BR}{{\mathbb {R}}}

      \newcommand{\BZ}{{\mathbb {Z}}}
      \newcommand{\CB}{{\mathcal {B}}}
     \newcommand{\CC}{{\mathcal {C}}} \renewcommand{\CD}{{\mathcal {D}}}
      
      \newcommand{\CH}{{\mathcal {H}}}
     \newcommand{\CI}{{\mathcal {I}}} 
      
      \newcommand{\CN}{{\mathcal {N}}}
     \newcommand{\CO}{{\mathcal {O}}} 
      \newcommand{\CR}{{\mathcal {R}}}

      \newcommand{\RN}{{\mathrm {N}}}

      \newcommand{\RT}{{\mathrm {T}}}

     \newcommand{\disc}{{\mathrm{disc}}}

     \newcommand{\End}{{\mathrm{End}}}  
\newcommand{\Frob}{{\mathrm{Frob}}}
     \newcommand{\Gal}{{\mathrm{Gal}}} \newcommand{\GL}{{\mathrm{GL}}}
      \renewcommand{\Im}{{\mathrm{Im}}}

\newcommand{\Jac}{{\mathrm{Jac}}}
     \newcommand{\Ker}{{\mathrm{Ker}}} 
      
     \newcommand{\fa}{{\mathfrak{a}}} 
      
  \newcommand{\n}{{\mathfrak{n}}} \newcommand{\m}{{\mathfrak{m}}}
\newcommand{\p}{{\mathfrak{p}}}  \newcommand{\h}{{\mathfrak{h}}}
\newcommand{\e}{{\mathfrak{e}}}\newcommand{\q}{{\mathfrak{q}}}
\newcommand{\E}{{\mathfrak{E}}} 
   \newcommand{\g}{{\mathfrak{g}}}
\renewcommand{\c}{{\mathfrak{c}}}
      
     \newcommand{\ord}{{\mathrm{ord}}}  

     \newcommand{\PGL}{{\mathrm{PGL}}} \newcommand{\Pic}{\mathrm{Pic}}
     \renewcommand{\Re}{{\mathrm{Re}}} 
     \newcommand{\Res}{{\mathrm{Res}}} 
     \newcommand{\RTr}{{\mathrm{Tr}}}  
       \newcommand{\fp}{{\mathfrak {p}}}

      \newcommand{\tr}{{\mathrm{tr}}}
     \newcommand{\tor}{{\mathrm{tor}}}

      \newcommand{\Supp}{{\mathrm{Supp}}}
     \font\cyr=wncyr10  \newcommand{\Sha}{\hbox{\cyr X}}
     \newcommand{\wt}{\widetilde} \newcommand{\wh}{\widehat}

     \newcommand{\sfrac}[2]{\left( \frac {#1}{#2}\right)}
     \newcommand{\nequiv}{\equiv\hspace{-9.5pt}/}
     
     \newcommand{\ov}{\overline}
     
     \newcommand{\lra}{\longrightarrow}
     
     \renewcommand{\lto}{\longmapsto}
     
     \theoremstyle{plain}
     \newtheorem{thm}{Theorem}[section] 
     \newtheorem{lem}[thm]{Lemma}  \newtheorem{prop}[thm]{Proposition}

     \renewcommand{\theequation}{\arabic{equation}}
   \numberwithin{equation}{section}

\begin{document}
\currannalsline{162}{2005} 

\title{Twisted Fermat curves over \\ totally real fields}

 \acknowledgements{}
\twoauthors{Adrian Diaconu}{Ye Tian}

 \institution{Columbia  
University, New York, NY\\
\email{cad@math.columbia.edu} 
\\ \vglue-9pt
McGill University, Montreal, Quebec, Canada\\
\email{tian@math.mcgill.ca}}


 \shorttitle{Twisted Fermat curves over  totally real fields}

\section{Introduction}

Let $p$ be a prime number, $F$ a totally real field such that
$[F(\mu_p): F]=2$ and $[F:\BQ]$ is odd. For $\delta \in F^\times$,
let $[\ \delta\ ]$ denote its class in $F^\times/F^{\times p}$. In
this paper, we show

\demo{\scshape Main Theorem} {\it There are
infinitely many classes $[\ \delta\ ]\in F^\times/F^{\times p}$
such that the twisted affine Fermat curves
$$W_\delta: \quad X^p+Y^p=\delta$$ have no $F$\/{\rm -}\/rational points.} 
 
\demo{Remark} It is clear that if $[\ \delta\ ]=[\ \delta'\ ]$,
then $W_\delta$ is isomorphic to $W_{\delta'}$ over $F$.  For any
$\delta \in F^\times,$ $W_\delta/F$ has rational points locally
everywhere.
\Enddemo

To obtain this result, consider the smooth open affine curve:
$$C_\delta: V^p=U (\delta -U),$$ and the morphism:
$$\psi_\delta: W_\delta \lra C_\delta;  \quad (x, y)\lto (x^p,
xy).$$ Let $C_\delta \rightarrow J_\delta$ be the Jacobian
embedding of $C_\delta/F$ defined by the point $(0, 0)$. We will
show that:

\begin{enumerate}
\item  If $L(1, J_\delta/F)\neq 0$, then $J_\delta (F)$ is a finite  
group
(cf.\ Theorem 2.1. of \S
2).

\quad The proof is based on Zhang's extension of the Gross-Zagier
formula to totally real fields and on Kolyvagin's technique of
Euler systems. One might use techniques of congruence of modular
forms to remove the restriction that the degree $[F:\BQ]$ is odd.
\item There are infinitely many classes $[\ \delta\ ]$ such that $L(1,
J_\delta/F)\neq 0$ (cf.\ Theorem 3.1. of \S 3; see also 2.2.4.).

\quad The proof is based on the theory of double Dirichlet series. The
condition that $[F(\mu_p):F]=2$ is essential for the technique we
use here.
\end{enumerate}
Combining $(1)$ and $(2)$,  one can see that the set
$$\Pi:=\Big\{[\ \delta\ ]\in F^\times/F^{\times p} \ \Big | \ J_\delta  
(F)\ \text{is torsion} \Big\}$$ is
infinite.

\Subsec{Proof of the Main Theorem assuming $(1)$ and $(2)$}
For any $\delta \in F^\times$, consider the twisting isomorphism
(defined  over $F(\sqrt[p]{\delta})$):
$$\iota_\delta: C_\delta \lra C_1; \quad (u, v)\lto (u/\delta,
v/\sqrt[p]{\delta^2}).$$ Define  $\eta_\delta: J_\delta \lra J_1$
to be the homomorphism associated to $\iota_\delta$.

Let $\Sigma_\delta$ denote the set $\iota_\delta \left(
C_\delta(F)\right)$. It is easy to see that:
\begin{enumerate}
\item[(i)] $\Sigma_\delta=\Sigma_{\delta'}$, if $[\ \delta\ ]=[\  
\delta'\ ]$,
\item [(ii)] $\Sigma_\delta \cap \Sigma_{\delta'}=\{(0, 0), (1, 0)\}$,
otherwise.
\end{enumerate}
For any $\delta \in F^\times$ with $[\ \delta\ ]\in \Pi$, and $[\
\delta\ ]\neq 1$, the diagram
$$\begin{array}{cccccccccccccc}
W_\delta (F)  \ \stackrel{\psi_\delta}{\longrightarrow} &C_\delta
(F)&
\hookrightarrow &J_\delta (F)\\
&\Biggr\downarrow\lefteqn{\iota_\delta} &&\Biggr  
\downarrow\lefteqn{\eta_\delta}\\
&C_1(F(\sqrt[p]{\delta}))&\hookrightarrow&
J_1(F(\sqrt[p]{\delta}))
\end{array}$$commutes.

Since the set
$$\bigcup_{\delta\in F^\times} J_1 (
F(\sqrt[p]{\delta}))_{\tor}\ \subset\  J_1(\ov{F})$$ is finite by the 
Northcott theorem, the set $\displaystyle{\bigcup_{[\ \delta\ ]\in
\Pi} \Sigma_\delta}$ is finite. Thus, for all but finitely many
$[\ \delta\ ]\in \Pi\setminus \{[1] \}$,\quad $\Sigma_\delta=\{(0,
0), (1, 0)\}$, and therefore $W_\delta$ has no $F$-rational
points. \hfill\qed

\demo{Remark} Our method is, in fact, effective: for any
$[\ \delta\ ]\in F^\times/F^{\times p}$, let
$$\Supp^{(p)}\left([\ \delta\ ]\right)=\left\{ \fp\ \text{prime of  
$F$}\ \Big|\ p\nmid v_\fp(\delta)\right\}.$$ Let
$L'$ be the Galois closure of $F(\mu_p)$, and let $S$ be the set
of places of $F$ above $2D_{L'/\BQ}$, where $D_{L'/\BQ}$ is the
discriminant of $L'/\BQ$. If $\Supp^{(p)}\left([\ \delta\
]\right)$ is not contained in $S$ and $L(1, J_\delta)\neq 0$, then
the twisted Fermat curve $W_\delta$ has no $F$-rational points
(see Proposition 2.2).

\demo{Acknowledgment}  We would like to thank D. Goldfeld, S.  
Friedberg,\break  J.
Hoffstein, H. Jacquet, V. A. Kolyvagin, L. Szpiro  for their help
and encouragement, and the referees for useful remarks and
suggestions. In particular, we are grateful to S. Zhang, who
suggested the problem to us,    for many helpful conversations.
The second author was partially supported by the Clay Mathematics
Institute.

  \section{\bf Arithmetic methods}

Fix $\delta \in F^\times \cap \CO_F$ such that $(\delta, p)=1$.
Let $\zeta=\zeta_p$ be a primitive $p$-th root of unity. The
abelian variety $J_\delta$ is absolutely simple, of dimension
$\displaystyle{g=\frac{p-1}{2}}$, and has complex multiplication
by $\BZ[\zeta]$ over the field $F(\mu_p)$. In this section we
show:
\begin{thm}
If $L(1, J_\delta/F)\neq 0${\rm ,} then $J_\delta (F)$ is finite.
\end{thm}

{\em Notation}. In this section,  for an abelian group $M$, set
$\wh{M}=M\otimes_\BZ \prod_p\BZ_p$ where $p$ runs over all primes.
For any ring $R$, let $R^\times$ denote the group of invertible
elements. For any ideal $\fa$ of $F,$ denote the norm
$\RN_{F/\BQ} (\fa)$  by $\RN \fa$. Let $\BA$ denote the adele ring
of $F$, and $\BA_f$ its finite part. Sometimes, we shall not
distinguish a finite place from its corresponding prime ideal.

\Subsec{The Hilbert newform associated to $J_\delta$} We first
recall some facts about $L$-functions of  twisted Fermat curves over
arbitrary number fields (see \cite{GR}, \cite{W2}). Let $F$  be any
number field, $L=F(\mu_p),$ $L_0=\BQ(\mu_p),$ and $F_0=L_0\cap F$.

For any place $w$ of $L$, denote by $w_0$ and  $v$   its  
restrictions to
$\BQ(\mu_p)$ and  $F$, respectively. Let $\chi_{w_0}$ and
$\chi_w$ be the $p$-th power residue symbols on $L_0^\times$ and
$L^\times,$ respectively, given by class field theory. Then
$\chi_w=\chi_{w_0}\circ \RN_{L/\BQ(\mu_p)}$.  The Jacobi sum
$$j(\chi_w, \chi_w)=-\sum_{\substack{{a\in \CO_L/w}\\{ a\neq 0, 1}}}  
\chi_w(a)\chi_w(1-a)$$ is an
integer in $L_0$ satisfying $j(\chi_w, \chi_w)=j(\chi_{w_0},
\chi_{w_0})^{i_{w/w_0}}$ and the Stickelberger relation:
$$\left(j(\chi_{w_0},
\chi_{w_0})\right)=\prod_{i=1}^{\frac{p-1}{2}}\sigma_i^{-1}(w_0)$$as
an ideal in $L_0.$ Here, $i_{w/w_0}$ is the inertial degree for
$w/w_0$, and $\sigma_i\in \Gal (L_0/\BQ)$ is the image of $i$
under the isomorphism $(\BZ/p\BZ)^\times \lra \Gal(L_0/\BQ)$.

Since $\delta \in \CO_F$ is coprime to $p$,  $C_\delta$ has good
reduction at $w$ for any $w\nmid p\delta$.   We know that the  
zeta-function
of the reduction $\wt{C_\delta}$ of $C_\delta$ at a place
$v$ of $F$ is
$$Z(\wt{C_\delta}, T)=\frac {P_v(T)}{(1-T)(1-\RN v T)},$$
with $$P_v(T)=\prod_{w|v}\prod_{\sigma} (1-\chi_w(\delta^2)^\sigma
j(\chi_w, \chi_w)^\sigma T^{f_v}),$$ where $f_v$ is the order of
$\RN v$ modulo $p,$ and $\sigma$ runs over representatives in
$\Gal(\BQ(\mu_p)/\BQ)$ of $\Gal(F_0/\BQ)$. Then the number of
points on $\tilde{J_\delta}$ (the reduction of $J_\delta$ at $v$)
is $P_v(1)$.

Now we give a bound on torsion points of $J_\delta(F)$. Let $F'$
be the Galois closure of $F/\BQ$, and assume that $F\cap
L_0=F'\cap L_0.$ This assumption is satisfied if $F$ is as in the
main theorem, or $F$ is Galois over $\BQ$. Let $L'=F'(\mu_p),$ and
let $q\nmid 2D_{L'/\BQ}$ be a prime. Let $\ell$ be a prime for
which there exists a place $w'| \ell$ of $L'$ such that
$\Frob_{L_0/F_0}(w'|_{L_0})$ is a generator of $\Gal(L_0/F_0)$,
$\Frob_{F'/F_0}(w'|_{F'})=1$ and
$\Frob_{\BQ(\mu_q)/\BQ}(w'|_{\BQ(\mu_q)})=1$. Then,
$\ell\equiv1\mod q$.  Let $v,$ $w$ and $w_0$ be the places of $F,$
$L$ and $L_0$, respectively, below $w'$. Then, $v$ is inert in
$L/F$ and $i_{w/w_0}=1.$ We have
$$P_v(1)=\prod_\sigma (1-\chi_w(\delta^2)^\sigma j(\chi_w,
\chi_w)^\sigma).$$ Since $v$ is inert in $L/F$ and $\delta\in
F^\times,$ we have $\chi_w(\delta^2)=1.$ Using the Stickelberger
relation and the fact that $j(\chi_{w_0}, \chi_{w_0})\equiv1\mod
(1-\zeta_p)^2$, one can show that\break\vskip-12pt\noindent  $j(\chi_w, \chi_w)=-\ell^f,$ for
$f=\frac{p-1}{2[F_0:\BQ]}.$ Then,
$P_v(1)=(1+\ell^f)^{[F_0:\BQ]}\equiv 2^{[F_0:\BQ]}\mod q.$
Consequently, there are no $q$-torsion points in $J_\delta (F)$.

Similarly, for the case  $q|2 D_{L'/\BQ}$, let $c_q\geq 1$ be the
smallest positive integer such that there is a $\sigma\in
\Gal(L'(\mu_{q^{c_q}})/\BQ)$ for which $\sigma|_L$ is a generator
of $\Gal(L/F)$, $\sigma|_{F'}=1$, and the restriction of $\sigma$
to $\Gal(\BQ(\mu_{q^{c_q}})/\BQ)$ has order greater than
$f=\frac{p-1}{2[F_0:\BQ]}.$ Then, $P_v(1)\nequiv\ 0\mod
q^{c_q[F_0:\BQ]}$. Let $M$ be defined\break\vskip-10pt\noindent by
$M:=\prod_{q|2D_{L'/\BQ}}q^{c_q[F_0:\BQ]}.$ It follows that
$J_\delta (F)_\tor \subset J_\delta [M],$ the subgroup of
$M$-torsion points of $J_\delta (\ov{F})$.

Let $F$ be a totally real field as in the main theorem. We have:
\begin{prop}Let $S$ be the set of places of $F$ above $2D_{L'/\BQ}$. If  
$\ \Supp^{(p)}\left([\ \delta\ ]\right)$
is not contained in $S$ and $L(1, J_\delta/F)\neq 0,$ then the
twisted\break Fermat curve $W_\delta$ has no $F$-rational points.
\end{prop}

Let $F$ be as in the introduction. Then $F_0=\BQ(\mu_p)^+$ is the
maximal totally real subfield of $L_0=\BQ(\mu_p).$ By the
reciprocity law, one can see that $w\mapsto \chi_w(\delta^2)$
defines a Hecke character, which we denote by $\chi_{[\delta^2]}.$
It depends only on\underline{} the class of $\delta^2$ and has
conductor above $\delta.$ By Weil \cite{W2}, the map $w\mapsto
j(\chi_w, \chi_w)\RN_{L/\BQ} w^{-\frac{1}{2}}$ also defines a
Hecke character on $L,$ denoted by $\psi,$ which has conductor
above $p.$ Thus, we have a (unitary) Hecke character on~$L$,
$$\chi_{[\delta^2]}\psi: \BA_L^\times \lra \BC^\times,$$
which is not of the form $\phi\circ \RN_{L/F},$ for any Hecke
character $\phi$ over $F.$ Then, there exists a unique holomorphic
Hilbert newform $f/F$ of pure weight $2$ with trivial central
character such that,
$$L_v(s, f/F)=\prod_{w|v}L_w(s-1/2, \chi_{[\delta^2]}\psi),$$ for all  
places $v$ of $F.$
Actually, the field over $\BQ$ generated by the Hecke eigenvalues
attached to $f$  is $F_0=\BQ(\mu_p)^+,$ and for the CM abelian
variety $J_\delta,$ we have
$$\begin{aligned}
L(s, J_\delta/F)&=\prod_{\sigma \in
\Gal(L_0/\BQ)\big/\Gal(L_0/F_0)} L(s-1/2, \chi_{[\delta^2]}^\sigma
\psi^\sigma)\\ \\&=\prod_{\sigma: F_0\hookrightarrow \BC} L(s,
f^\sigma/F).
\end{aligned}$$Note that $L(s, J_\delta)$ only depends on the
class $[\ \delta\ ]$ of $\delta,$ and the above equality
holds for any local factor.

\Subsec{A nonvanishing result} Let $\pi$ be the automorphic
representation associated to $f,$ and let $N$ be its conductor.
Let $S_0$ be any finite set of places of $F,$ including all
infinite places and the places dividing $N.$ Choose a quadratic
Hecke character $\xi$ corresponding to a totally imaginary
quadratic extension of $F,$ unramified at $N,$ where $\xi
(N)\cdot(-1)^g=-1$ (since $F$ is of odd
degree, we have $(-1)^g=-1$); i.e., the epsilon factor of $L(s,
\pi\otimes \xi)$ is $-1.$ Let $\CD(\xi; S_0)$ denote the set of
quadratic characters $\chi$ of $F^\times/\BA_F^\times,$ for which
$\chi_v=\xi_v,$ for all $v\in S_0.$ With the above notation and
assumptions, by a theorem of Friedberg and Hoffstein \cite{FH}, there
exist infinitely many quadratic characters $\chi\in \CD (\xi;
S_0)$ such that $L(s, \pi\otimes \chi)$ has a simple zero at the
center $s = 1/2.$

Choose such a $\chi,$ and let $K$ be the totally imaginary
quadratic extension of $F$ associated to it. The conductor of
$\chi$ is coprime to $N,$ and the $L$-function  $L(s, f/K)=L(s-1/2,
\pi)L(s-1/2, \pi\otimes \chi)$ has a simple zero at $s=1.$ Let $d$  
denote
the discriminant of $K/F.$

\Subsec{Zhang\/{\rm '}\/s formula}
\vglue-8pt
\Subsubsec{The $(N, K)$-type Shimura curves}
Let $\CO$ be the subalgebra of $\BC$ over $\BZ$ generated by the
eigenvalues  of $f$ under the Hecke operators. In
our case, $\CO=\BZ[\zeta+\zeta^{-1}]$ is the ring of integers of $F_0.$
In \cite{Zhang1} (see also \cite{Ca1},  \cite{Ca2}),
Zhang constructs a Shimura curve $X$ of $(N, K)$-type, and proves
that there exists a unique abelian subvariety $A$ of the Jacobian
$\Jac (X)$ of dimension $[\CO:\BZ]=g,$ such that
$$L_v(s, A)=\prod_{\sigma: \CO \hookrightarrow \BC} L_v(s,
f^\sigma/F),$$for all places $v$ of $F.$ By the construction of
$f,$ it follows that $L_v(s, A/F)=L_v(s, J_\delta/F)$ for all
places $v$ of $F.$ Therefore, by the isogeny conjecture proved by
Faltings, $A$ is isogenous to $J_\delta$ over $F.$ In particular,
the complex multiplication by $\CO\subset \BQ(\mu_p)^+$ on $A$ is
defined over $F.$

Now, let us recall the constructions of $X$ and $A.$

The $L$-function of $\pi\otimes \chi$ satisfies the functional
equation
$$L(1-s, \pi\otimes \chi)=(-1)^{\left|\Sigma\right|} \RN_{F/\BQ}  
(Nd)^{2s-1}
L(s, \pi\otimes \chi),$$where $\Sigma=\Sigma (N, K)$ is the
following set of places of $F:$
$$\Sigma (N, K)=\left\{ v\ \Big |\  v|\infty,\ \text{or}\
\chi_v(N)=-1\right\}.$$ Since the sign of the functional equation
is $-1,$ by our choice of $K,$ the cardinality of $\Sigma$ is odd.
Let $\tau$ be any real place of $F.$ Then, we have:
\begin{enumerate}
\item  Up to isomorphism, there exists a unique quaternion algebra $B$
  such that $B$ is ramified at exactly the places in $\Sigma \backslash
\{\tau\}$;
\item There exist embeddings $\rho: K \hookrightarrow B$ over $F.$
\end{enumerate}

 From now on, we fix an embedding $\rho: K\rightarrow B$ over $F.$

Let $G$ denote the algebraic group over $F,$ which is an inner
form of $\PGL_2$ with $G(F)\cong B^\times/F^\times.$ The group
$G(F_\tau)\cong \PGL_2(\BR)$ acts on $\CH^{\pm}=\BC\setminus \BR.$
Now,  for any open  compact subgroup $U$ of $G(\BA_f),$ we have an
analytic space
$$ S_{U}(\BC)=G(F)_+ \backslash \CH^+ \times G(\BA_f)/U,$$
where $G(F)_+$ denotes the subgroup of elements in $G(F)$ with
positive determinant via $\tau.$

Shimura has shown that $S_U(\BC)$ is the set of complex points of
an algebraic curve $S_U,$ which descends canonically to $F$ (as a
subfield of $\BC$ via $\tau$). The curve $S_U$ over $F$ is
independent of the choice of $\tau.$

There exists  an order $R_0$ of $B$ containing $\CO_K$ with reduced
discriminant~$N.$ One can choose $R_0$ as follows. Let $\CO_B$ be
a maximal order of $B$ containing $\CO_K,$ and let $\CN$ be an ideal of
$\CO_K$ such that $$\RN_{K/F}\CN \cdot \disc_{B/F}=N,$$ where
$\disc_{B/F}$ is the reduced discriminant of $\CO_B$ over $\CO_F.$
Then, we take $$R_0=\CO_K+\CN \cdot \CO_B.$$ Take $U=\prod_v
R^\times_v/\CO_v^\times.$ The corresponding Shimura curve $X:=S_U$
is compact.

Let $\xi\in \Pic (X)\otimes \BQ$ be the unique class whose degree
is $1$ on each connected component and such that,
$$\RT_m\xi=\deg (\RT_m)\xi,$$for all integral ideals $m$ of $\CO_F$
coprime to $Nd.$ Here, the $\RT_m$ are the Hecke operators.

\Subsubsec{Gross-Zagier-Zhang formula}
Now, we define the basic class in $\Jac (X)(K)\otimes \BQ,$ where
$\Jac (X)$ is the connected component of $\Pic (X),$ from the
CM-points on the curve $X.$ The CM points corresponding to $K$ on
$X$ form a set:
$$\CC: \ G(F)_+\setminus G(F)_+\cdot h_0\times G(\BA_f)/U\cong
T(F)\setminus G(\BA_f)/U; \qquad [(h_0, g)]\leftrightarrow
[g],$$where $h_0\in \CH^+$ is the unique fixed point of the torus
$T(F)=K^\times/F^\times.$ \medbreak

For a CM point $z=[g]\in \CC,$ represented by $g\in G(\BA_f),$ let
$$\Phi_g: K\lra \wh{B}, \qquad t\lto g^{-1}\rho (t) g.$$Then, $\End
(z):=\Phi_g^{-1}(\wh{R_0})$ is an order of $K,$ say $\CO_n=\CO_F+n
\CO_K,$ for a (unique) ideal $n$ of $F.$ The ideal $n,$ called the
conductor of $z,$ is independent of the choice of the representative
$g.$ By Shimura's theory, every CM point of conductor $n$ is
defined over the abelian extension $H_n'$  of $K$ corresponding to
$K^\times\setminus \wh{K}^\times/ \wh{F}^\times \wh{\CO}_n^\times$
via class field theory.

Let $P_1$ be a CM point in $X$ of conductor $1,$ which is defined
over $H_1',$ the abelian extension of $K$ corresponding to
$K^\times \setminus \wh{K}^\times /\wh{F}^\times
\wh{\CO}_K^\times.$ The divisor $P= \Gal (H_1'/K)\cdot P_1$
together with the Hodge class defines a class
$$x:=[P-\deg (P) \xi]\in \Jac (X)(K)\otimes \BQ,$$where $\deg P$
is the multi-degree of $P$ on the geometric components. Let $x_f$
be the $f$-typical component of $x.$ In \cite{Zhang2},
Zhang generalized the Gross-Zagier formula to the totally real field case, by  
proving that

$$L'(1,f/K)=\frac{2^{g+1}}{\sqrt{\RN(d)}}\cdot \|f\|^2\cdot \|x_f\|^2,$$
where $\|f\|^2$ is computed on the invariant measure on
$$\PGL_2(F)\setminus \CH^g\times \PGL_2(\BA_f)/U_0(N)$$induced by  
$dxdy/y^2$ on $\CH^g,$ and where
$$U_0(N)=\left\{ \begin{pmatrix}a&b\\c&d\end{pmatrix} \in
\GL_2(\wh{\CO}_F)\big | c\in \wh{N}\right\}\subset
\GL_2(\wh{F}),$$ and $\|x_f\|^2$ is the Neron-Tate pairing of
$x_f$ with itself.

\Subsubsec{The equivalence of nonvanishing of  $L$-factors}
For any \hbox{$\sigma: F\hookrightarrow \BC,$}  it is known by a result of  
Shimura
that $L(1, f/F)\neq 0$ is equivalent to $L(1, f^\sigma /F)\break\neq 0.$
One can also show this using Zhang's formula above. To see this,
assume $L(1, f/F)\neq 0.$ Then, $\|x_f\|\neq 0,$ and therefore,
$\|x_{f^\sigma}\|\neq 0.$ It follows that $L'(1, f^\sigma /K)\neq
0.$ Since $L(1, f/F)\neq 0,$   the $L$-function $L(s,
f^\sigma/F)$ has a positive sign in its functional equation. Thus,
$L(1, f^\sigma/F)\neq 0.$ In fact, to obtain our main theorem, we do
not need this equivalence, but we may see that Theorem 3.1 is
equivalent to statement (2) in the introduction.

\Subsec{The Euler system of CM points}
We now assume that $L(1, \chi_{[\delta^2]}\psi)\neq 0,$ or
equivalently,  $L(1, f/F)\neq 0.$ Then by the equivalence of
nonvanishing of $L(1, f^\sigma)$ for all embeddings $\sigma:
F\hookrightarrow \BC,$ we have that $L(1, J_\delta/F)\neq 0.$ By
Zhang's formula, we also know that $\|x_f\|\neq 0.$

Let $\CN$ be the set of square-free integral ideals of $F$ whose
prime divisors are inert in $K$ and coprime to $Nd.$ For any $n\in
\CN,$ define
$$H_n=\prod_{\ell |n} H'_\ell\subset H_n', \qquad H_1=H_1'.$$ Let $u_n$  
denote the cardinality
of $(\widehat{\CO}_n^\times\cap K^\times
\widehat{F}^\times)/\widehat{\CO}_F^\times.$ Then, $H_\ell/H_1$ is
a cyclic extension of degree
$t(\ell)=\frac{\RN(\ell)+1}{u_1/u_\ell}.$
\vskip2pt

For each $n\in \CN,$ let $P_n$ be a CM point of order $n$ such
that $P_n$ is contained in $\RT_\ell P_m$ if $n=m\ell \in \CN$ and
$\ell$ is a prime ideal of $F.$ Let $y_n=\RTr_{H_n'/H_n}\pi
(P_n)\in A(H_n),$ where $\pi$ is a morphism from $X$ to $\Jac
(X)$ defined by a multiple of the Hodge class.

The points  $\{y_n\}_{n\in \CN}$ form an Euler system (see
\cite[Prop.~7.5]{Tian}, or \cite[Lemma 7.2.2]{Zhang1}) so that, for
any $n=m\ell \in \CN$ with $\ell$ a prime ideal of $F,$ 
\begin{enumerate}
\item $\displaystyle{{u_n}^{-1} \sum_{\sigma\in \Gal(H_n/H_m)}
y_n^\sigma={u_m}^{-1}a_\ell y_m}$;
\item For any prime ideal $\lambda_m$ of $H_m$ above $\ell,$ and for 
$\lambda_n$  the unique prime above $\lambda_m,$
$$\Frob_{\lambda_m} y_m\equiv  y_n \mod \lambda_n;$$
\item  The class $x_f$ is equal to $y_K:=\tr_{H_1/K}y_1$ in  
$\big(A(K)\otimes \BQ\big)\big/\BQ^\times.$
\end{enumerate}

Theorem 2.1 follows with the nontrivial Euler system by
Kolyvagin's standard argument (see \cite{K2}, \cite{KL2},
\cite{G}, and \cite[Th.~A]{Zhang1}).

\vglue-16pt
\phantom{up}

\section{\bf Analytic methods}

Let $r=4$ or an odd prime, and let $L=F(\zeta_r),$ with
$[L:F]=2.$ Let $\psi$ be a unitary Hecke character of $L.$ In this
section, we show:

\begin{thm}
There are infinitely many classes $\delta \in F^\times/F^{\times
r}$ such that $L\left(\frac{1}{2}, \chi_{[\ \delta\ ]}\psi\right)$
does not vanish.
\end{thm}

Let $\rho$ be a unitary Hecke character of $F.$ The purpose of
this section is to construct a perfect double Dirichlet series
$Z(s, w; \psi; \rho)$ similar to an Asai-Flicker-Patterson type
Rankin-Selberg convolution, which possesses meromorphic
continuation to $\BC^2$ and functional equations. Then, Theorem
$3.1$ will follow from the analytic properties of $Z(s, w; \psi;
\rho)$ (when $r=4$, see [7]). To do this, it is necessary to recall the Fisher-Friedberg
symbol in \cite{FF1}.

\Subsec{The $r$\/{\rm -th}\/ power residue symbol}
Let $S'$ be a finite set of non-archimedean places of $L$
containing all places dividing $r,$ and such that the ring of
$S'$-integers $\CO_L ^{S'}$ has class number one. We shall also
assume that $S'$ is closed under conjugation and that $\psi$ and
$\rho$ are both unramified outside $S'.$

Let $S_\infty$ denote the set of all archimedean places of $L,$
and set $S = S'\cup S_\infty.$ Let $I_{L}(S)$ (resp. $\CI_{L}(S)$)
denote the group of fractional ideals (resp. the set of all
integral ideals) of $\CO_L$ coprime to $S'.$ In \cite{FF1}, Fisher
and Friedberg have shown that the $r$-th order symbol $\chi_n$
can be extended to $I_{L}(S)$ i.e., $\chi_\frak{n}(\frak{m})$ is
defined for $\frak{m},$ $\frak{n}\in I_{L}(S).$ Let us recall
their construction.

For a non-archimedean place $v \in S',$ let $\frak{P}_v$ denote
the corresponding ideal of $L.$ Define $\frak{c}=\prod_{v\in S'}
\frak{P}_v^{r_v}$ with $r_v=1$ if $\ord_v(r)=0$, and $r_v$
sufficiently large such that, for $a\in L_v,$ $\ord_v(a-1)\geq r_v$
implies that $a\in (L^\times_v)^r$. Let $P_L(\c)\subset I_L(S)$ be
the subgroup of principal ideals $(\alpha)$ with $\alpha \equiv 1
\mod \c,$ and let $H_\frak{c}=I_L(S)/P_L(\c)$ be the ray class
group modulo $\c.$ Set $R_\c=H_\c \otimes \BZ/r\BZ,$ and write the
finite group $R_\c$ as a direct product of cyclic groups. Choose a
generator for each, and let $\E_0$ be a set of ideals of $\CO_L,$
prime to $S,$ which represent these generators. For each $\e_0\in
\E_0,$ choose $m_{\e_0} \in L^\times$ such that
$\e_0\CO_L^{S'}=m_{\e_0}\CO_L^{S'}.$ Let $\E$ be a full set of
representatives for $R_\c$ of the form $\prod_{\e_0\in \E_0}
\e_0^{\lambda_{\e_0}}.$ Note that $\e\CO_L^{S'}=m_\e \CO_L^{S'}$
for all $\e\in \E.$ Without loss, we suppose that $\CO_L^{S'}\in
\E$ and $m_{\CO_L^{S'}}=1.$

Let $\m, \n\in I_L(S)$ be coprime. Write $\m=(m)\e\g^r$ with $\e
\in \E,$ $m\in L^\times$, $m\equiv 1\mod \c$ and $\g\in I_L(S),$
$(\g, \n)=1.$ Then the $r$-th power residue symbol $\sfrac{m m_\e
}{\n}_r$ is defined. If $\m=(m')\e'\g^{'r}$ is another such
decomposition, then $\e'=\e$ and $\sfrac{m' m_{\e'}
}{\n}_r=\sfrac{m m_\e }{\n}_r.$

In view of this, the $r$-th power residue symbol
$\sfrac{\m}{\n}_r$ is defined to be\break\vskip-12pt\noindent $\sfrac{m m_\e }{\n}_r,$ and
the character $\chi_\m$ is defined by $\chi_\m
(\n)=\sfrac{\m}{\n}_r.$ This extension of the $r$-th power
residue symbol depends on the above choices. Let $S_\m$ denote the
support of the conductor of $\chi_\m.$ It can be easily checked
that if $\m=\m'\fa^r$, then $\chi_\m(\n)=\chi_{\m'}(\n)$ whenever
both are defined. This allows one to extend $\chi_\m$ to a
character of all ideals of $I_L(S\cup S_\m).$

The extended symbol possesses a reciprocity law: if $\m , \n\in
I_L(S)$ are coprime, then $\alpha (\m, \n)=\chi_\m (\n)\chi_\n
(\m)^{-1}$ depends only on the images of $\m, \n$ in~$R_\c.$

In our situation, we also need the following lemma:
\begin{lem}
The natural morphism
$$I_F(S)/P_F(\c)\lra I_L(S)/P_L(\c)$$
has kernel of order a power of $2.$
\end{lem}

\Proof 
If $[\n]$ is in the kernel, i.e., $\n=(\alpha)$ in $I_L(S)$ is a
principal ideal with $\alpha \equiv 1 \mod \c,$ then $\alpha
/\ov{\alpha}$ is a root of unity with $\alpha/\ov{\alpha} \equiv 1
\mod \c.$ Now let $W$ be the set of roots of unity in $L$ which
are $\equiv 1\mod \c.$ Let $W_0$ be the subset of $W$ of elements
of the form $u/\ov{u}$ for some unit $u$ in $\CO_L$ and $u\equiv
1\mod \c.$ It is clear that $W_0\supset W^2.$ Then, the map
$$\Ker \left( I_F(S)/P_F(\c)\rightarrow I_L(S)/P_L(\c) \right)
\lra W/W_0; \qquad \n \lto \alpha/\ov{\alpha}$$ is obviously
injective; i.e., the order of the kernel of the natural map in
this lemma is  a power of $2.$
\Endproof\vskip4pt

Since $r$ is odd, using the lemma, we may choose a suitable set
$\E_0$ of representatives since the beginning such that if  $\m\in
I_F(S),$ then the decomposition $\m=(m)\e\g^r$ is such that $m\in
F^\times$, $\e, \g\in I_F(S).$

Using the symbol $\chi_\n,$ we shall construct a perfect double
Dirichlet series $Z(s, w; \psi; \rho)$ (i.e., possessing
meromorphic continuation to ${\BC}^2$) of type:
$$Z(s, w; \psi; \rho)\;=\;Z_{S}(s, w; \psi; \rho)\;\;=  
*\sum_{\frak{n}\in
\CI_{F}(S)} L_{S}(s, \psi\,
\chi_{\frak{n}})\,\rho(\frak{n})\,\RN_{F/\BQ}(\frak{n})
^{-w},\leqno (3.1)$$ where the sum is over the set of all integral
ideals of $\CO_F$ coprime to $S',$ for $\frak{n}\in \CI_{F}(S)$
square-free, the function $L_{S}(s, \psi\, \chi_{\frak{n}})$ is
precisely the Hecke $L$-function attached to $\psi\,
\chi_{\frak{n}}$ with the Euler factors at all places in $S$
removed, and where $*$ is a certain normalizing factor. For an arbitrary
$\frak{n}\in \CI_F(S),$ write $\frak{n} = \frak{n}_1\frak{n}_2
^{r}$ with $\frak{n}_1$ $r$-th power free. If $L_{S}(s,
\psi\,\chi_{\frak{n}_1})$ denotes the Hecke $L$-series associated
to $\psi\,\chi_{\frak{n}_1}$ with the Euler factors at all places
in $S$ removed, then $L_{S}(s, \psi\,\chi_{\frak{n}})$ is defined
as $L_{S}(s, \psi\,\chi_{\frak{n}_1})$ multiplied by a Dirichlet
polynomial whose complexity grows with the divisibility of
$\frak{n}$ by powers (see $(3.10),$ $(3.12)$ and $(3.13)$ for
precise definitions).

Based on the analytic properties of $Z(s, w; \psi; \rho)$, we show
the following result which is stronger than Theorem 3.1.

\begin{thm}  {\rm 1)}  There exist infinitely many $r$\/{\rm -th}\/ power free ideals  
$\n_1$ in $\CI_F(S)$ with trivial image
in $R_\c$ for which the special value $L_S(\frac{1}{2}, \chi_\n
\psi)$ does not vanish.
\smallbreak
{\rm 2)}  Let $\kappa_\c$ denote the number of  characters of
$R_\frak{c}$ whose restrictions to $F$ are also  characters of the
ideal class group of $F${\rm ,} and let $\kappa$ be the residue of the
Dedekind zeta function $\zeta_F(s)$ at $s=1$. Then for
$x\rightarrow \infty,$  
$$\sum_{\substack{{\RN_{F/\BQ}(\frak{n}) < x}\\{\frak{n}\in
\CI_F(S)}\\{\n\,=\, (n)}\\{[\n]\,=\,1}}}L_{S}\left(\frac{1}{2},\,
\chi_{\frak{n}}\psi\right)\;\sim \; \frac{\kappa\cdot \kappa_\c}{h_F
\cdot|R_{\frak{c}}|}\,\frac{L_{S}(1, \psi)\,L_{S}( \frac{r}{2},
\psi^{r})}{L_S(\frac{r}{2}+1, \psi^r)
}\prod_{\substack{{v\;\text{in}\;F}\\{v\in S'}}}\left(1 - q_v
^{-1}\right)\cdot x, \ \leqno (3.2)$$ where $[\n]$ denotes the
image of the ideal $\n$ in $R_\frak{c}$.
\end{thm}

  {\em Remarks.} i)  By the above definition of the extended  
$r$-th
power residue symbol, it is easy to see that the first part of
this theorem is equivalent to Theorem $3.1.$

\smallbreak
ii)  In fact, by a well-known result of Waldspurger \cite{Wa}, it
will follow that $L_S(\frac{1}{2}, \chi_\n \psi)\geq 0,$ for $\n
\in \CI_F(S),\, \n=(n)$ and trivial image in $R_\frak{c}.$ We will
see this in the course of the proof of   Theorem $3.3.$
\smallbreak
iii)  Following \cite{DGH}, by a simple sieving process, one can
prove the more familiar variant of the above asymptotic formula
where the sum is restricted to square-free principal
ideals.

\Subsec{The series $Z_{\rm aux}(s, w; \psi; \rho)$
and metaplectic Eisenstein series}
To obtain the correct definition of $Z(s, w; \psi; \rho)$,  let
$G_{0}(\n,\, \m),$ for $\m,$ $\n\in \CI_{L}(S),$ be given by
$$G_{0}(\n,\, \m)\;\;\;=\prod_{\substack{ v\\
{\ord_v(\n)=k} \\ {\ord_v(\m) = l}}} G_{0}(\frak{p}_v^k,
\frak{p}_v^l),\ \leqno (3.3)$$ where, for $k,$ $l\geq 0,$
$$G_{0}(\p_v^k,\, \p_v^l)\,=\,\begin{cases} 1& \text{if $l=0,$}\\
q_v^{\frac{k}{2}}& \text{if $k + 1=l;$ $l\not\equiv 0\pmod r,$}\\
-\,q_v^{\frac{k - 1}{2}}& \text{if $k + 1 = l;$ $l > 0;$ $l \equiv  
0\pmod r,$}\\
q_v ^{\frac{l}{2} - 1}(q_v - 1)& \text{if $k\geq l;$ $l > 0;$ $l \equiv  
0\pmod r,$}\\
0 & \text{otherwise.}
\end{cases}\ \leqno (3.4)$$ Here $q_v$ denotes the absolute value of the  
norm of
$v.$ Also, let $G(\chi_{\m_{1}} ^{*})$ (where $\m_{1}$ denotes the
$r$-th power free part of $\m$ and $\chi_{\frak{a}}
^{*}(\frak{b}): = \chi_{\frak{b}}(\frak{a})$)  be the normalized
Gauss sum appearing in the functional equation of the (primitive)
Hecke $L$-function associated to $\chi_{\m} ^{*}.$ If $\n ^{*}$
denotes the part of $\n$ coprime to $\m_{1},$ then set $$G(\n,
\m)\,:=\,\overline{\chi_{\m_{1}} ^{*}(\n^*)}\,G(\chi_{\m_{1}}
^{*})\,G_{0}(\n, \m).$$ Now, let $\psi$ be as above. For $\n\in
\CI_{L}(S)$ and $\Re(s) > 1,$ let $\Psi_{S}(s, \n, \psi)$ be the
absolutely convergent Dirichlet series defined by
$$\Psi_{S}(s, \n, \psi)\;=\; L_S \left(rs - \frac{r}{2} +
1, \psi ^r \right)\sum_{\m\in \CI_{L}(S)}\frac{\psi(\m) G(\n,
\m)}{\RN_{L/\Bbb Q}(\m) ^s}.$$ This series can be realized as a
Fourier coefficient of a metaplectic Eisenstein series on the
$r$-fold cover of $\GL(2)$ (see \cite{KP} and \cite{K}). It
follows as in Selberg \cite{S}, or alternatively, from Langlands'
general theory of Eisenstein series \cite{L} that $\Psi_{S}(s, \n,
\psi)$ has meromorphic continuation to $\BC$ with only one
possible (simple) pole at $s = \frac{1}{2} + \frac{1}{r}.$
Moreover, this function is bounded when $|\Im(s)|$ is large in
vertical strips, and satisfies a functional equation as $s\to 1 -
s$ (see Kazhdan-Patterson  \cite[Cor.\ II.2.4]{KP}).

For $\Re(s),\, \Re(w)
> 1,$ let $Z_{\rm aux}(s, w; \psi; \rho)$
be the auxiliary double Dirichlet series defined by
$$Z_{\rm aux}(s, w; \psi; \rho)\;\;
=  \sum_{\frak{n}\in \CI_F(S)} \frac{\Psi_{S}(s, \n,
\psi)\rho(\frak{n})}{\RN_{F/\BQ}(\frak{n}) ^w}.\leqno (3.5)$$ Let
$\tilde{\rho}$ be the Hecke character of $L$ given by
$\tilde{\rho} = \rho\,\circ\, \RN_{L/F}.$ As we shall shortly see,
$Z_{\rm aux}(s, w; \psi\,\tilde{\rho}; \ov{\rho})$ is the type
of object that constitutes a building block in the process of
constructing the perfect double Dirichlet series $Z(s, w; \psi;
\rho).$ Set
$$\Gamma_{\rm aux} ^{*}(s, \psi\, \tilde{\rho}) \;= \prod_{v\in
S_\infty}\;  \prod_{j = 1} ^{r - 1} L_v \Big(s - \frac{1}{2} +
\frac{j}{r}, \psi_{v}\, \tilde{\rho}_{v} \Big),$$ and let
$$\widehat{Z}_{\rm aux}(s, w; \psi\,\tilde{\rho}; \bar{\rho})\,
:=\,\Gamma_{\rm aux} ^{*}(s, \psi\, \tilde{\rho})\cdot
Z_{\rm aux}(s, w; \psi\,\tilde{\rho}; \bar{\rho}).$$ Let
$\CR_{1}$ be the tube region in $\BC ^{2}$ whose base $\CB_{1}$ is
the convex region in $\BR ^{2}$ which lies strictly above the
polygonal contour determined by $(0,2),$ $(1,1)$, and the rays
$y=-2x+2$ for $x\le 0$ and $y=1$ for $x\ge 1$. As a simple
consequence of the analytic properties of $\Psi_{S}(s, \n, \psi)$
($\n\in \CI_{L}(S)$), we have the following:

\begin{prop} The double Dirichlet series
$Z_{\rm aux}(s, w; \psi\,\tilde{\rho}, \bar{\rho})$ is
holomorphic in $\CR_{1},$ unless $\psi ^{r}\tilde{\rho} ^{r} = 1$
when it has only one simple pole at $s = \frac{1}{2} +
\frac{1}{r}$. Furthermore{\rm ,} $\widehat{Z}_{\rm aux}(s, w;
\psi\,\tilde{\rho}, \bar{\rho})$ satisfies the functional equation
$$\begin{aligned} &\widehat{Z}_{\rm aux}(s, w; \psi\,\tilde{\rho},  
\bar{\rho})
\,\cdot\prod_{v\in S'} \left( 1 - (\psi \tilde{\rho})
^{-r}(\pi_v)\,q_v ^{rs - \frac{r}{2} - 1}\right)\\
&\hskip70pt =\; \sum_{\eta,\, \tau}\,A_{\eta,\, \tau} ^{(\psi,\,
\rho)}(1 - s)\, \widehat{Z}_{\rm aux}(1 - s, 2s + w - 1; \psi
^{-1}\tilde{\rho} ^{-1}\eta, \psi\,\rho\,\tau),
\end{aligned}\leqno (3.6)$$ where each $A_{\eta,\, \tau} ^{(\psi,\,
\rho)}(s)$ is a polynomial in the variables $q_v ^{s},\, q_v
^{-s}$ $(v\in S'),$ and the sum is over a finite set of id\'ele
class characters $\eta$ and $\tau,$ unramified outside $S$ and with 
orders dividing $r.$
\end{prop}
\vglue-16pt \pagegoal=50pc

\Subsec{The double Dirichlet series $\wt{Z} (s, w; \psi;
\rho)$}
It turns out that the function $Z_{\rm aux}(s, w;
\psi\,\tilde{\rho}, \bar{\rho})$ possesses another functional
equation. To describe it, we introduce a new double Dirichlet
series $\wt{Z}(s, w; \psi; \rho)$ defined for $\Re(s),\, \Re(w)\break
> 1$ by  \setcounter{equation}{6}
\begin{eqnarray} &&\\
&&\wt{Z}(s, w; \psi; \rho)\; =\;L_S(rs + rw
+ 1 - r, \psi ^{r} \tilde{\rho} ^{r})\sum_{\substack{{\frak{m}\in
\CI_{L}(S)}\\{\m-\text{imaginary}}}} \frac{\psi(\m)\,L_{S}(w,\,
\chi_{\m} ^{*}\,\rho)}{\RN_{L/\BQ}(\m) ^{s}}\nonumber\\&&\hskip 12pt\cdot
\sum_{\h\in \CI_F(S)} \frac{(\psi \rho)(\h)\,\chi_{\m}
^{*}(\h_{1})}{\RN_{F/\BQ}(\h) ^{2s - 1}\,\RN_{F/\BQ}(\h)
^{w}}\;\prod_{\substack{{v}\\ { \ord_{v}(\h_0)> 0}}}
\left[\,(\chi_{\m} ^{*}\, \rho)(\pi_{v})\,q_{v} ^{-w }\, -\; q_{v}
^{-1}\, \right]\nonumber\\
&&\hskip
12pt\cdot  \prod_{\substack{{v}\\ {\ord_{v}(\RN_{L/F}(\m)) >
0}\\{ \ord_{v}(\h_{2}) > 0}}} (1 - q_v ^{-1}) 
\prod_{\substack{{v-\text{split in}\ L}\\
{\ord_{v}(\RN_{L/F}(\m)) = 0}\\ {\ord_{v}(\h_{2}) > 0}}}
\left[\,(\chi_{\m} ^{*}\, \rho)(\pi_{v})\,q_{v} ^{- w - 1}\, +\,
1\, - \, 2 q_v ^{-1}\,\right]\nonumber\\*
&&\hskip
12pt\cdot\prod_{\substack{{v-\text{inert in}\
L}\\{ \ord_{v}(\h_2) > 0}}} \left[\,1\,-\;(\chi_{\m} ^{*}\,
\rho)(\pi_{v})\,q_{v} ^{- w - 1}\,\right].\nonumber
\end{eqnarray}
 In the above formula,  an ideal $\frak{m}\in  
\CI_{L}(S)$ is called
{\em imaginary}, if it has no divisor in $\CI_F(S),$ other than
$\CO_F$. The function $L_{S}(w,\, \chi_{\m} ^{*}\,\rho)$
represents the $L$-series defined over $F$ (not necessarily
primitive) associated to $\chi_{\m} ^{*}\,\rho$ with the Euler
factors corresponding to places removed in $S$. Also, all the
products are over places of $F$, $\pi_{v}$ is the local parameter
of $F_{v}$ ($F_{v}$ denoting the completion of $F$ at $v$), and
$q_v$ is the absolute value of the norm in $F$ of $v.$

Let $\CR_{2}$ denote the tube region in $\BC ^{2}$ whose base
$\CB_{2}$ is the convex region in $\BR ^2$ which lies strictly
above the polygonal contour determined by $(1,1),$
$(\frac{3}{2},0)$ and the rays $y=-x+\frac{3}{2}$ for $y\le 0$ and
$x=1$ for $y\ge 1.$ Recall that $L_{S}(w,\, \chi_{\m} ^{*}\,\rho)$
differs from a primitive $L$-series by only finitely many Euler
factors (i.e., the factors corresponding to places in $S$ and to
places $v$ for which $\ord_v(\RN_{L/F}(\m))\equiv 0 \pmod r$).
Applying the functional equation of $L_{S}(w,\, \chi_{\m}
^{*}\,\rho)$ and some standard estimates, one can easily show that
the function $\widetilde{Z}(s, w; \psi; \rho)$ is holomorphic in
$\CR_{2},$ unless $\rho = 1$ where it has only one simple pole at
$w = 1$. The following proposition gives the functional equation
connecting the double Dirichlet series $Z_{\rm aux}(s, w;
\psi\,\tilde{\rho}, \bar{\rho})$ and $\widetilde{Z}(s, w; \psi;
\rho).$
\setcounter{equation}{7}

\begin{prop} The function $\wt{Z}(s, w; \psi;
\rho)$ is holomorphic in $\CR_{2},$ unless $\rho$ is the trivial
character when it has a simple pole at $w = 1$. Furthermore{\rm ,} for
$\Re(s),\, \Re(w) > 1,$ there exist the functional equations
\begin{multline}
\prod_{v\in S_\infty} L_v \left(1 - w, \rho_{v}
\right)\,\cdot\prod_{v\in S'} \left(1 - \rho ^{-r}(\pi_v)\,q_v
^{-rw}\right)\cdot \widetilde{Z}(s + w - {\scriptstyle
\frac{1}{2}}, 1 - w; \psi; \rho)\\  = \prod_{v\in
S_\infty} L_v \left(w, \rho_{v} ^{-1}
\right)\,\cdot\sum_{\tau}B_{\tau} ^{(\rho)}(w)\, Z_{\rm aux}(s,
w; \psi \tilde{\rho}\, \tau, \bar{\rho}),
\end{multline}
 and $$\begin{aligned} &\prod_{v\in S_\infty}  
L_v
\left(w, \rho_{v} ^{-1} \right)\,\cdot\prod_{v\in S'} \left(1 -
\rho ^{r}(\pi_v)\,q_v ^{rw - r}\right)\cdot Z_{\rm aux}(s, w;
\psi \tilde{\rho}, \bar{\rho})\\& \hskip 28pt = \prod_{v\in
S_\infty} L_v \bigl(1 - w, \rho_{v}
\bigr)\,\cdot\sum_{\tau}C_{\tau} ^{(\rho)}(1 - w)\,
\widetilde{Z}(s + w - {\scriptstyle\frac{1}{2}}, 1 - w;
\psi\,\tau; \rho),
\end{aligned}\leqno (3.9)$$
where{\rm ,} as before, $B_{\tau} ^{(\rho)}(w),$ $C_{\tau} ^{(\rho)}(w)$
are polynomials in the variables $q_v ^{w},\, q_v ^{-w}$ $(v\in
S').$ The above products are over the places of $k$ corresponding
to those in $S,$ and the sums are over a finite set of id{\rm \'{\hskip-6pt\it e}}le
class characters $\tau,$ unramified outside $S$ and orders
dividing $r.$
\end{prop}
 
The proof of this proposition will be given in the next section.

Let $\alpha$ and $\beta$ be the involutions on $\BC ^2$ given by
$$\alpha : (s, w) \rightarrow (1 - s, 2s + w - 1) \quad \hbox{and}
\quad \beta : (s, w) \rightarrow (s + w - {\scriptstyle
\frac{1}{2}}, 1 - w).$$ It can be easily checked that these
involutions generate the dihedral group $D_{8}$ of order $8.$ It
follows directly from Propositions $3.2$ and $3.3$ that both\break
$\widetilde{Z}(s + w - {\scriptstyle\frac{1}{2}}, 1 - w; \psi;
\rho)$ and $Z_{\rm aux}(s, w; \psi \tilde{\rho}, \bar{\rho})$
can be continued to $\CR_{1}\,\cup\,\CR_{2}$.  Clearly, this
applies to $Z_{\rm aux}(s, w; \psi, \rho)$ (replace $\psi$ by
$\psi \tilde{\rho} ^{-1}$ and $\rho$ by $\bar{\rho}$). It follows
from the functional equation $(3.6)$ that $Z_{\rm aux}(s, w;
\psi \tilde{\rho}, \bar{\rho})$ can be continued to
$\CR_{1}\,\cup\,\CR_{2}\,\cup\,\alpha(\CR_{2}),$ and hence, by
$(3.8),$ the function $\widetilde{Z}(s + w - {\scriptstyle
\frac{1}{2}}, 1 - w; \psi; \rho)$ continues to this region. The
double Dirichlet series $Z_{\rm aux}(s, w; \psi\wt{\rho},
\ov{\rho})$ may have only one simple pole in $\CR_2$, namely
$w=1$, and this pole occurs only if $\rho$ is the trivial
character. This fact follows easily by inspection of the proof of
Proposition 3.3 (see \S 3.1). Then from the functional equation
(3.6), one can see that $Z_{\rm aux}(s, w; \psi\wt{\rho},
\ov{\rho})$ may have a pole only at $w=2-2s$ in $\alpha (\CR_2)$,
provided $\psi^r|_{\CO_F} \cdot \rho^r$ is trivial. The last fact
also applies to $\wt{Z}(s+w-\frac{1}{2}, 1-w; \psi, \rho)$, by the
functional equation $\beta$ in (3.8).

\Subsec{The double Dirichlet series $ Z(s, w; \psi;
\rho)$} To define the perfect double Dirichlet series $Z(s, w;
\psi; \rho),$ let $L_{S}(s, \chi_{\frak{n}}\psi),$ for $\n\in
\CI_{F}(S),$ be given by $$L_{S}(s,\, \chi_{\frak{n}}\psi)\,:=\,
L_{S}(s,\, \chi_{\n_{1}} \psi)P_{\n}(s,\, \psi),$$ where $\n_1$
denotes the $r$-th power free part of $\n,$ and $P_{\n}(s, \psi)$
is the Dirichlet polynomial defined by \setcounter{equation}{9}
\begin{eqnarray}
&&\\
&&  P_{\n}(s,\, \psi) =  
 \prod_{\substack{{v}\\{\ord_v(\n_{1})
> 0}}}\Biggr(1 + \psi(\pi_{v})
\,q_{v} ^{1 - 2s} + \cdots + \psi(\pi_{v}) ^{\ord_v(\n) - 1} q_{v}
^{(\ord_v(\n) - 1)(1 - 2s)}\Biggr)\nonumber\\
&&\cdot \prod_{\substack{{v}\\{\ord_v(\n)= r \mu}\\{v-\text{inert
in}\;L}}}\Biggr( \Big(1 - \psi (\pi_{v})\, q_{v} ^{- 2
s}\Big)\left(1 + \psi(\pi_{v})\,q_{v} ^{1 - 2s} + \cdots\right.\nonumber\\
&&\left.\qquad\qquad +
\psi(\pi_{v}) ^{r\mu -
1}\,q_{v} ^{(r\mu - 1)(1 - 2s)}\right) + \psi(\pi_{v}) ^{r\mu}\,q_{v} ^{r\mu (1 -
2s)}\left(1 + q_{v} ^{- 1}\right)\Biggr)\nonumber \\
&&\cdot \prod_{\substack{{v}\\{\ord_v(\n)= r \omega}\\{v = v'
\bar{v}'\;\text{in}\;L}}}\Biggr( (1-(\chi_{\n_{1}}
\psi)(\pi_{v'})\, q_{v} ^{-s})(1 -
(\chi_{\n_{1}}\psi)(\pi_{\bar{v}'})\,q_{v}
^{-s})\big (1 + \psi(\pi_{v})\,q_{v} ^{1 - 2s} + \cdots\nonumber\\
&&\qquad\qquad + \psi(\pi_{v}) ^{r\omega - 1}\,q_{v} ^{(r\omega -
1)(1 - 2s)}\big) + \psi(\pi_{v}) ^{r\omega}\,q_{v} ^{r\omega (1 -
2s)}\left(1 - q_{v} ^{- 1}\right)\Biggr).\nonumber
\end{eqnarray}
 Here the products are over places $v$ of  
$F$, and $\pi_{v}$
denotes the local parameter of $F_{v}.$ It can be seen that
these polynomials satisfy a functional equation as $s\to 1 - s,$
and that we have the estimate $$P_{\n}(s,\, \psi)\ll_{\varepsilon}
\RN_{F/\BQ}(\n) ^{\varepsilon}\;\;\;\;\;\;\;\;\;\;(\varepsilon
> 0,\; \Re(s)\ge {\scriptstyle \frac{1}{2}} ). \leqno (3.11)$$  
  Furthermore,
if $\psi(\overline{\m})=\overline{\psi(\m)},$ for $\m \in
\CI_{L}(S),$ then $P_{\n}(s,\, \psi) \geq 0,$ for $s\in \BR$.
Later, we shall specialize $\psi$ to be (essentially) a
normalized Jacobi sum, which obviously satisfies this property.

For $\Re(s),\, \Re(w) > 1,$ we define $Z(s, w; \psi; \rho)$  as
\setcounter{equation}{11}
\begin{eqnarray} \qquad
Z(s, w; \psi; \rho)&=&Z_S(s, w; \psi; \rho)\\
&=&L_S(rs + rw + 1 - r,\, \psi ^{r} \tilde{\rho}
^{r})\sum_{\frak{n}\in \CI_F(S)} \frac{L_{S}(s,\,
\chi_{\frak{n}}\psi)\rho(\frak{n})}{\RN_{F/\BQ}(\frak{n})
^w}.\nonumber
\end{eqnarray}
 Applying the functional equation and the
convexity bound of $L_{S}(s, \chi_{\frak{n}}\psi)$\break ($\n\in
\CI_{F}(S)$), we see that $Z(s, w; \psi; \rho)$ is holomorphic
in $\CR_{1},$ if the character $\psi ^r$ is nontrivial.
Representing the normalizing factor $L_S(rs + rw + 1 - r, \psi
^{r} \tilde{\rho} ^{r})$ by its Dirichlet series, then after
multiplying and reorganizing, we can write $Z(s, w; \psi; \rho)$
as
$$Z(s, w; \psi; \rho)\;\;=\sum_{\frak{n}\in
\CI_{F}(S)} \frac{L_{S}(s, \,\chi_{\n_{1}}\psi)\,Q_{\n}(s,\,
\psi)\,\rho(\frak{n})}{\RN_{F/\BQ}(\frak{n}) ^w}, \leqno (3.13)$$
where $Q_{\n}(s, \psi),$ for $\n\in \CI_F(S),$ is a new set of
Dirichlet polynomials which can be easily expressed in terms of
$P_{\n}(s, \psi).$\vglue2pt

Referring to the definition of $\widetilde{Z}(s, w; \psi; \rho)$
given in $(3.7)$, replace $L_{S}(w, \chi_{\m} ^{*}\,\rho)$ by its
Dirichlet series, the sum being over $\n,$ say. For fixed $\m \in
\CI_L(S)$ imaginary, and $\n \in \CI_F(S),$ collect the terms
contributing to $(\chi_{\m} ^{*}\,\rho)(\n)\,\RN_{F/\BQ}(\n)
^{-w}.$ Switching the order of summation, we obtain:

\begin{prop} For $\Re(s),\, \Re(w) > 1${\rm ,} 
$$Z(s, w; \psi; \rho)\,=\,L_{S}(2s, \psi)\widetilde{Z}(s, w; \psi;
\rho),\leqno (3.14)$$ where the $L$-function is defined over $F$.
\end{prop}

Assuming both $\psi ^r$ and $\psi ^{r}\tilde{\rho} ^{r}$ to be
nontrivial, we see  from Proposition $3.4$ that
$$L_{S}(2s + 2w - 1, \psi) \widetilde{Z}(s + w -
{\scriptstyle\frac{1}{2}}, 1 - w; \psi; \rho)$$ continues to
$\beta(\CR_{1}),$ and hence, from the above discussion, it
continues to $\CR_{1}\,\cup
\beta(\CR_{1})\,\cup\,\CR_{2}\,\cup\,\alpha(\CR_{2})$.  Note that
the convex closure of this tube region is $\BC ^{2}$. As $\psi
^{r}\tilde{\rho} ^{r} \neq 1$, and therefore, by Propositions
$3.2$ and $3.3$, the function $\widetilde{Z}(s + w - {\scriptstyle
\frac{1}{2}}, 1 - w; \psi; \rho)$ does not have a pole at $s =
\frac{1}{2} + \frac{1}{r}$, one can easily check that the only
possible poles of $L_{S}(2s + 2w - 1, \psi) \widetilde{Z}(s + w -
{\scriptstyle \frac{1}{2}}, 1 - w; \psi; \rho)$ are the
hyperplanes $w = 0$ and $w=2- 2s$. Clearly, both are  simple poles,
and they may occur only if $\rho$ and $\psi^r|_{\CO_F}\cdot
\rho^r$    are both  trivial.

Consequently, by the convexity theorem for holomorphic functions
of several complex variables (see \cite{Ho}) and by Proposition
$3.4,$ we have the following:

\begin{thm} When $\psi ^r$ and $\psi ^{r}\tilde{\rho}
^{r}$ are nontrivial{\rm ,} the function
$$(w - 1)(2s + w - 2)Z(s, w; \psi; \rho)$$ has analytic
continuation to $\BC^{2},$ and for any fixed $s,$ it is \/{\rm (}\/as a
function of the variable $w${\rm )} of order one.
\end{thm}

The fact that, for any fixed $s,$ the above function is of order
one follows as in \cite[Prop.~$3.11$]{DGH}.

By Proposition $3.4$ and $(3.7),$ one finds that, for $\Re(s)
> \frac{1}{2},$
\vglue12pt\noindent (3.15)
\begin{small}
\begin{eqnarray*}
&&\underset{w=1}{\Res}\  Z(s, w; \psi; 1)  = 
L_{S}(2s, \psi)\,L_S(rs + 1,
\psi ^{r})\nonumber \\
&&\qquad \cdot \prod_{\substack{{v\;\text{in}\;F}\\{v\in
S'}}}\Biggr[\left(1 - q_v ^{-1}\right)
\sum_{\substack{{\frak{m}\in \CI_{L}(S)}\\{\m-\text{imaginary}}}}
\Biggr( \frac{\kappa\,\psi(\m) ^{r}\,\prod_{v'|\m}\left(1 - q_{v'}
^{-1}\right)}{\RN_{L/\BQ}(\m)
^{rs}}  \sum_{\h\in \CI_F(S)} \frac{\psi (\h)
^r}{\RN_{F/\BQ}(\h) ^{2rs}}\nonumber\\
&&\qquad \cdot \prod_{\substack{{
v}\\{\ord_{v}(\RN_{L/F}(\m))
> 0}\\{\ord_{v}(\h)
> 0}}}
  \big (\,1\, -\, q_v ^{-1}\,\big ) \prod_{\substack{{ v-\text{split
in}\ L}\\{\ord_{v}(\RN_{L/F}(\m)) = 0}\\{\ord_{v}(\h) > 0}}} \big
(\,1\, - \, q_v ^{-1}\,\big ) ^2\prod_{\substack{{ v-\text{inert
in $L$}}\\{\ord_{v}(\h) > 0}}} \big ( 1\,-\;q_{v} ^{-
2}\,\big)\Biggr)\Biggr]\nonumber \\&&\quad =\kappa
L_{S}(2s, \psi)\,L_S(rs, \psi
^{r})\prod_{\substack{{v\;\text{in}\;F}\\{v\in S'}}}\left(1 - q_v
^{-1}\right),
\nonumber
\end{eqnarray*}
\end{small}

\noindent 
 where $\kappa$ denotes the residue at $w =  
1$ of the
Dedekind zeta-function $\zeta_F(w).$

We are now in the position to give the proof of Theorem $3.3.$

\demo{Proof of Theorem $3.3$} As before, let $\rho = \prod
\rho_v$ be a unitary Hecke character of $F$ unramified outside
$S.$ We further assume that $\rho$ is of finite order. For
$\Re(s),\, \Re(w) > 1,$ consider the double Dirichlet series
$Z_{1}(s, w; \psi; \rho)$ defined by $$Z_{1}(s, w; \psi;
\rho)\;\;\,=\; \sum_{\substack{{\n \in \CI_F(S)}\\{\n\,=\,
(n)}\\{[\n]\,=\,1}}} \frac{L_{S}(s,
\,\chi_{\n_{1}}\psi)\,Q_{\n}(s,\,
\psi)\,\rho(\frak{n})}{\RN_{F/\BQ}(\n) ^w}. \leqno (3.16).$$
By expressing this function as $$Z_{1}(s, w; \psi; \rho)\;=\;
\frac{1}{h_F\cdot|R_{\frak{c}}|}\,\sum_{\rho_{1},\,\rho_{2}}\,Z(s,
w; \psi; \rho\rho_{1}\wh{\rho}_{2}),$$ where $\rho_1$ ranges over
the characters of the ideal class group of $F$, $\rho_2$ ranges
over the characters of $R_{\frak{c}},$ and $\wh{\rho}_{2}$ is the
restriction of $\rho_2$ to $F$,  it follows from Theorem $3.5$
that $Z_{1}(s, w; \psi; \rho)$ is holomorphic on $\BC ^2,$ except
for $w = 1$ and $w=2 - 2s$, where it might have simple poles.
Furthermore,
$$\lim_{w\rightarrow 1}(w - 1) ^2\, Z_{1}({\scriptstyle \frac{1}{2}},  
w; \psi;
\rho)\;\;\;=\lim_{(s, w) \rightarrow (\frac{1}{2}, 1)}(w - 1)(2s +
w - 2)Z_{1}(s, w; \psi; \rho)\,=\,0,$$ and, therefore,
$Z_{1}(\frac{1}{2}, w; \psi; 1)$ has at most a simple pole at $w =
1.$ To compute its residue, recall the functional equation
satisfied by $L(s, \chi_{\n_{1}}\psi)$ with $\n_{1}\in \CI_F(S)$
$r$-th power free (see \cite[Ch.\ VII, \S 7]{W1}). Combining
this with the functional equation of the polynomial $Q_{\n}(s,
\psi)$ ($\n\in \CI_F(S)$), we find that
\begin{eqnarray*}
 L_{S}(s, \,\chi_{\n_{1}}\psi)\,Q_{\n}(s,\,
\psi) &=&\varepsilon(s,\, \chi_{\n_{1}}\psi)\cdot L_{S}(1 - s,
\,\chi_{\n_{1}}\psi)\,Q_{\n}(1 - s,\, \psi)\\
&&\cdot\prod_{v\in
S_{\infty}}\frac{L_{v}(1 - s,\, \psi_{v})}{L_{v}(s,\,
\psi_{v})}\,\cdot\prod_{v\in S'}\frac{L_{v}\left(1 - s,\,
(\chi_{\n_{1}}\psi)_{v}\right)}{L_{v}\left(s,\,
(\chi_{\n_{1}}\psi)_{v}\right)}.\end{eqnarray*}
 A simple local
computation shows that $\varepsilon(\frac{1}{2}, \chi_{\n_{1}}\psi)
=\psi(\n) \varepsilon(\frac{1}{2}, \psi)$. It immediately follows
that $Z_{1}(s, w; \psi; 1)$ satisfies the functional equation
$$\begin{aligned} &\prod_{v\in S_\infty} L_v (s, \psi_{v})\,
\cdot\prod_{v\in S'} \Big (1 - \psi ^{r}(\pi_v)\,q_v ^{rs - r}\Big
)\cdot Z_{1}(s, w; \psi; 1)\\& \qquad =\; \prod_{v\in S_\infty}
L_v (1 - s, \psi_{v})\,\cdot\sum_{\rho}D_{\rho} ^{(\psi)}(1 - s)\,
Z_{1}(1 - s, 2s + w - 1; \psi; \rho),
\end{aligned} \leqno (3.17)$$ where $D_{\rho} ^{(\psi)}(s)$
are polynomials in the variables $q_v ^{s},\, q_v ^{-s},$ $v\in
S',$ and the sum is over a finite set of id\'ele class characters
$\rho,$ unramified outside $S$ and orders dividing $r$. As $r$ is
odd, and $\psi$, restricted to the group of principal ideals of $F$,
is quadratic and nontrivial, it follows that $Z_1(s, w; \psi; 1)$
does not have a pole at $w=2-2s$. Then $(3.15)$ yields
$$\underset{w=1}{\Res}\  Z_1\left(\frac{1}{2}, w; \psi;  
1\right)=\frac{\kappa \cdot
\kappa_\c}{h_F \cdot |R_\frak{c}|} L_S(1, \psi) L_S \left(\frac{r}{2},
\psi^r\right) \prod_{\substack{{v \, \text{in}\, F}\\{v\in
S'}}}(1-q_v^{-1}),\leqno (3.18)$$where $\kappa_\c$ denotes the
number of  characters of $R_\frak{c}$ whose restrictions to $F$
are also characters of the ideal class group of $F$.

To complete the proof, we define the double Dirichlet series
$Z_0(s, w; \psi; \rho)$ by simply replacing in $(3.16)$ the
polynomial $Q_\n (s, \psi)$ by  $P_\n (s, \psi)$ defined in
$(3.10)$. Note that
$$Z_0(s, w; \psi;\rho)=\frac{1}{h_F \cdot
|R_\frak{c}|}\sum_{\rho_1, \rho_2} \frac{Z(s, w; \psi;
\rho\rho_1\rho_2)}{L_S(rs+rw+1-r, \psi^r \wt{\rho}^r
\wt{\rho}^r_1)},$$and therefore, $Z_0(s, w; \psi;\rho)$ may have
additional poles at the zeros of the incomplete $L$-functions
$L_S(rs+rw+1-r, \psi^r \wt{\rho}^r \wt{\rho}^r_1)$. It is
well-known that these zeros occur in the region $\Re(s + w) < 1.$ In
particular, the function $Z_0(\frac{1}{2}, w; \psi; 1)$ is
holomorphic for $\Re(w)>\frac{1}{2},$ except for $w=1$, where it
has a simple pole. Using $(3.18)$, we can compute its residue as
$$\underset{w=1}{\Res}\  Z_0\left( \frac{1}{2}, w; \psi;
1\right)\,=\,\frac{\kappa\cdot
\kappa_\c}{h_{F}\cdot|R_{\frak{c}}|}\, \frac{L_{S}(1,
\psi)\,L_{S}(\frac{r}{2}, \psi ^{r})}{L_S(\frac{r}{2}+1,
\psi^r)}\prod_{\substack{{v\;\text{in}\;F}\\{v\in S'}}}\left(1 -
q_v ^{-1}\right)>0.\leqno (3.19)$$

This implies that $L_S(\frac{1}{2}, \chi_{\n_1}\psi)\neq 0$ for
infinitely many $r$-th power free ideals $\n_1$ in $\CI_F(S)$ with
trivial image in $R_\c$, which is the first assertion of  
Theorem~3.3.

For the remaining part, one needs to apply a Tauberian theorem. To
keep the argument as simple as possible, note first that, as
$\psi(\ov{\m})=\ov{\psi(\m)},$ for $\m \in \CI_L(S),$ we have
$P_\n(s, \psi)\geq 0,$ for $s\in \BR.$ On the other hand, by the
comment made right after   Lemma $3.2,$ any $r$-th power free
ideal $\n_1$ in $\CI_F(S)$ with trivial image in $R_\c$ can be
decomposed as $\n_1=(n_1)\g^r$ with $n_1\in F^\times,\,\,
n_1\equiv 1\mod \c$ and $\g \in I_F(S).$ By definition, the
character $\chi_{\n_1}$ coincides with the classical $r$-th power
residue symbol $\chi_{n_1}$ given by class field theory. It
follows that the incomplete $L$-series $L_S(s, \chi_{\n_1}\psi)$
differs from the complete Hecke $L$-series associated to $L(s,
\chi_{n_1}\psi)$ by only finitely many local factors. Recall that
the latter is the $L$-series associated to a Hilbert modular
form. As the set $S'$ is closed under conjugation, it follows from
a well-known result of Waldspurger \cite{W1} that
$L_S(\frac{1}{2}, \chi_\n \psi)\geq 0,$ for $\n\in \CI_F(S),
\,\,\n=(n)$ and trivial image in $R_\c.$ Hence, the function
$Z_0({\scriptstyle \frac{1}{2}}, w; \psi; 1),$ for $\frak{R}
(w)>1,$ is given by a Dirichlet series with nonnegative
coefficients. The second part of   Theorem $3.3$ now follows
from the  Wiener-Ikehara Tauberian theorem.
\Endproof \vskip4pt

{\em Remark.} With some additional effort, one can exhibit an
error term on the order of $O (x^\theta)$ with $\theta <1$ in the
asymptotic formula $(3.2).$ Also, the remark following Theorem
$3.3$ implies that the Hecke $L$-series $L_S(\frac{1}{2},
\chi_{\n_1}\psi)\neq 0$ for infinitely many square-free principal
ideals $(n)$ in $\CI_F(S)$ with trivial image in~$R_\c.$ Any such
ideal has a generator $n\in F$ with $n\equiv 1 \mod
\c.$

\Subsec{Proof of Proposition $3.3$}
Recall that for $\frak{a} \in \CI_L(S),$ we defined
$\chi_{\frak{a}} ^{*}$ by $\chi_{\frak{a}} ^{*}(\frak{b}): =
\chi_{\frak{b}}(\frak{a})$ ($\frak{b} \in \CI_L(S)$). Note that
every ideal $\m$ of $\CO_L$ can be uniquely decomposed as $\m =
\m'\h,$ where $\m'$ is an imaginary ideal of $\CO_L,$ and $\h$ is
a real ideal; that is, $\h \in \CO_F.$ For $\m \in \CI_{L}(S)$
imaginary and $r$-th power free, let $\varepsilon(w, (\chi_{\m}
^{*}\,\rho) ^{-1})$ denote the epsilon-factor in the functional
equation of $L(w, (\chi_{\m} ^{*}\,\rho) ^{-1})$ (as a Hecke
$L$-function of $F$). Also, for $\m \in \CI_{L}(S)$ imaginary and
$\h \in \CI_F(S),$ coprime and $r$-th power free, let
$G(\chi_{\m\h} ^{*})$ be the normalized Gauss sum in the
functional equation of the Hecke $L$-function (of the field $L$)
associated to $\chi_{\m \h} ^{*},$ i.e., $\varepsilon(\frac{1}{2},
\chi_{\m\h} ^{*}).$ We set $\m_{0}$ and $\h_{0}$ to be the product
of all distinct prime ideals dividing $\m$ and $\h,$ respectively.

The following lemma is a consequence of a standard local
computation. The details will be omitted.

\begin{lem} Let $\m$ and $\h$ be integral ideals as above. Assume that
the images of $\m\h$ and $\m$ in $R_{\frak{c}}$ are $\frak{e}$ and
$\e',$ respectively. Then{\rm ,}
\begin{multline*}
G(\chi_{\m\h} ^{*})\,\varepsilon\left( \frac{1}{2},
(\chi_{\m} ^{*}\,\rho) ^{-1}\right) \\* =
\,C_{\e,\,\e',\,\rho}\cdot\,\eta(\e)
^{-1}\eta(\m_{1}\h_{1})\,\tilde{\rho} (\m_{0}) ^{-1}\,\chi_{\m}
^{*}(\h_{0})\,\chi_{\h} ^{*}(\m_{0})\,\chi_{\m}
^{*}(\overline{\m}_{0}) ^{-1},
\end{multline*}
where $\tilde{\rho} =
\rho\,\circ\, \RN_{L/F},$ $C_{\e,\,\e',\,\rho}$ is a constant
depending on just $\e,$ $\e'$ and $\rho,$ and $\eta$ is a Hecke
character unramified outside $S$ and order dividing $r.$
Furthermore{\rm ,} if $\e'$ is replaced by $\e''$ with $\e'/\e''$ a real
ideal{\rm ,} then both $C_{\e,\,\e',\,\rho}$ and $\eta$ do not change.
\end{lem}

\setcounter{equation}{19}
{\it Proof of Proposition $3.3$}.
Using $(3.5),$ we have
\begin{eqnarray}
&&\\
&& Z_{\rm aux}(s, w; \psi\,\tilde{\rho},
\bar{\rho})\nonumber\\
&&\qquad =\sum_{\frak{n}\in \CI_F(S)} \frac{\Psi_{S}(s, \n,
\psi\,
\tilde{\rho})\,\overline{\rho(\frak{n})}}{\RN_{F/\BQ}(\frak{n})
^w}\nonumber\\&&\qquad =  L_S\left(rs - \frac{r}{2} + 1, \psi ^{r} \tilde{\rho}
^{r}\right) \sum_{\substack{{\frak{m}\in \CI_{L}(S)}\\{\n \in
\CI_F(S) }}} \frac{(\psi\,
\tilde{\rho})(\m)\,\overline{\rho(\frak{n})}\,G(\n,
\m)}{\RN_{L/\BQ}(\m) ^s\,\RN_{F/\BQ}(\frak{n}) ^w}
\nonumber\\&&\qquad = L_S \left(rs -
\frac{r}{2} + 1, \psi ^{r} \tilde{\rho} ^{r}\right)
\sum_{\substack{{\frak{m}\in \CI_{L}(S)}\\{\n \in \CI_F(S) }}}
\frac{(\psi\,
\tilde{\rho})(\m)\,\overline{\rho(\frak{n})}\,\,\overline{\chi_{\m_{1}}
^{*}(\n^*)}\,G(\chi_{\m_{1}} ^{*})\,G_{0}(\n, \m)}{\RN_{L/\BQ}(\m)
^s\,\RN_{F/\BQ}(\n) ^w},\nonumber
\end{eqnarray}
 where $\n ^{*}$ denotes the part of $\n$  
coprime to
$\m_{1}.$ In the last sum, replace $\m$ by $\m\h$ with $\m \in
\CI_{L}(S)$ imaginary and $\h$ real. As we shall see, the only
contribution to the sum comes from $\m$ and $\h$ for which their
$r$-th power free parts $\m_1$ and $\h_1$ are coprime. Then, we
have
$$\begin{aligned} &\sum_{\substack{{\frak{m}\in \CI_{L}(S)}\\{\n \in
\CI_F(S) }}} \frac{(\psi\,
\tilde{\rho})(\m)\,\overline{\rho(\frak{n})}\,\,\overline{\chi_{\m_{1}}
^{*}(\n^*)}\,G(\chi_{\m_{1}} ^{*})\,G_{0}(\n, \m)}{\RN_{L/\BQ}(\m)
^s\,\RN_{F/\BQ}(\n) ^w}\;\;\;\;\,=\sum_{\substack{{\m \in
\CI_{L}(S)}\\{\m-\text{imaginary}}}} \frac{(\psi\,
\tilde{\rho})(\m)}{\RN_{L/\BQ}(\m) ^s}\\&\qquad\qquad\cdot
\sum_{\substack{{\h\in \CI_L(S)}\\{\n \in
\CI_F(S)}\\{\h-\text{real}}}} \frac{(\psi\,
\tilde{\rho})(\h)\,\overline{\rho(\frak{n})}\,\,\overline{\chi_{\m_{1}\h 
_{1}}
^{*}(\n^*)}\,G(\chi_{\m_{1} \h_{1}} ^{*})\,G_{0}(\n, \m
\h)}{\RN_{L/\BQ}(\h) ^s\,\RN_{F/\BQ}(\n) ^w}.
\end{aligned}\leqno (3.21)$$

Next, we separate the contribution of $\h$ in the inner sum. To do
so, let $\m_1$ denote the $r$-th power free part of an ideal
$\frak{m}\in \CI_{L}(S),$ and set  $\m_{0}$ to be the product of
all distinct prime ideals dividing $\m_{1},$ and
$$\m_{2}\;\;\;\;:=\prod_{\substack{{v}\\{ \ord_v(\m)=re_{v}}}}
\p_v^{re_{v}}.$$ For fixed $\m,$ $\n$ and $\h$ as above, let
$\p_v$ be a prime ideal of $L$ dividing $\h_{0}.$ Upon replacing
this prime ideal by its conjugate, we can assume that $\ord _v(\m)
= 0.$ Recall that
$$G_{0}(\n,\, \m)\;\;\;=\prod_{\substack{{ v}\\{
\ord_v(\n)=k}\\{\ord_v(\m)=l}}} G_{0}(\p_v^k, \p_v^l),$$ where
$G_{0}(\p_v^k,\, \p_v^l)$ is given by $(3.4).$ As $\ord_v(\m\h) =
\ord_v(\h)\not\equiv 0\pmod r$ (this condition implying that $
\ord_v(\n) = \ord_v(\h) - 1$), and $\n \in \CI_F(S),$ we can
decompose $\n = (\h/\h_{0}\h_{2})\n'$ with $\n'\in \CI_F(S)$
coprime to $\h_1.$ Also, we have
\begin{eqnarray*}
\ord_{v}(\n)&=&\ord_{\bar v}(\n)\,\ge\, \ord_{\bar v}(\m\h) - 1\\[4pt]
&=&
\, \ord_{\bar v}(\m) + \ord_{v}(\h) - 1\,=\, \ord_{\bar v}(\m) +
\ord_v(\n),
\end{eqnarray*}
 which implies $\ord_{\bar v}(\m) = 0.$ It
immediately follows that $\m$ and $\h_1$ are coprime. Then, by
$(3.4),$ we can write
\setcounter{equation}{21}
\begin{eqnarray}
 G(\chi_{\m_{1} \h_{1}}
^{*})\,G_{0}(\n, \m \h)\,&=&G(\chi_{\m_{1}\h_{1}}
^{*})\,G_{0}\left(\frac{\h}{\h_{0}\h_{2}},
\frac{\h}{\h_{2}}\right)G_{0}(\n', \m
\h_{2})\\
&=&G(\chi_{\m_{1}\h_{1}}
^{*})\,\RN_{L/\BQ}\left(\frac{\h}{\h_{0}\h_{2}}\right)
^{\frac{1}{2}}G_{0}(\n', \m \h_{2}).\nonumber
\end{eqnarray}
 Furthermore, we have
\begin{eqnarray*}
&& \!{\bf }{\bf } G_{0}(\n', \m \h_{2})  =\prod_{\substack{{ v}\\{
\ord_v(\n')=k_v}\\{\ord_v(\m)=l_v}\\{\ord_v(\h_{2})=re_v}}}
G_{0}(\p_v^{k_v}, \p_v ^{l_v + re_v})\\[4pt]
&&\! =\prod_{\substack{{
v}\\{l_v\not\equiv 0\,(r)}\\{k_v + 1= l_v + re_v}}}
G_{0}(\p_v^{k_v}, \p_v ^{l_v + re_v}) \cdot \prod_{\substack{{
v}\\{l_v\equiv 0\,(r)}\\{k_v + 1\ge l_v + re_v}}}
G_{0}(\p_v^{k_v}, \p_v ^{l_v + re_v})\\[4pt]
&&\! =\prod_{\substack{{
v}\\{l_v\not\equiv 0\,(r)}\\{k_v + 1 = l_v + re_v}}}
q_v^{\frac{(l_v - 1) + re_v}{2}} \cdot\prod_{\substack{{v}\\{l_v \equiv
0\,(r)}\\{k_v + 1 = l_v + re_v > 0}}} -\;\,q_v^{\frac{l_v + re_v -
2}{2}} \cdot\prod_{\substack{{ v}\\{l_v\equiv 0\,(r)}\\{k_v
\ge l_v + re_v
> 0}}} q_v ^{\frac{l_v + re_v}{2}}(1 - q_v
^{-1})\\[4pt]
&&\! =\RN_{L/\BQ}\left(\frac{\m\h_{2}}{\m_{0}}\right)
^{\frac{1}{2}} \cdot\prod_{\substack{{ v}\\{l_v \equiv 0\,(r)}\\{k_v + 1
= l_v + re_v
> 0}}} -\;\,q_v^{- 1}  \cdot\prod_{\substack{{v}\\{l_v\equiv 0\,(r)}\\{k_v \ge  
  l_v
+ re_v > 0}}} (1 - q_v ^{-1}).
\end{eqnarray*}
 One can decompose
$\n'$ as 
\begin{eqnarray*}
\n'& =& \n_{1}\cdot\,
\RN_{L/F}\left(\frac{\m}{\m_{0}}\right)\cdot\,\h_{2}\\
&&\cdot
\prod_{\substack{{ v-\text{complex}}\\{l_v\equiv 0\,(r);\;l_{\bar
v} = 0}\\{l_v + re_v > 0}\\{\alpha_{v} := 1 + k_v - l_v - re_v \ge
0}}} \RN_{L/F}(\p_v) ^{\alpha_{v} - 1}\;\;\;\;
\cdot\prod_{\substack{{ v-\text{real}}\\{e_v
> 0}\\{\beta_{v} := 1 + k_v - re_v \ge 0}}} \q_{v} ^{\beta_{v} - 1},
\end{eqnarray*}
 with
$\n_{1}$ coprime to $\m\h.$ Here, if $v$ is complex such that $l_v
= l_{\bar v} = 0,$ then one chooses either $v$ or $\bar v,$ but
not both. As $\n = (\h/\h_{0}\h_{2})\n',$ we also have 
\begin{eqnarray*}
\n &=&
\n_{1}\cdot\,
\RN_{L/F}\left(\frac{\m}{\m_{0}}\right)\cdot\,\frac{\h}{\h_{0}}\\*
&& \cdot
\prod_{\substack{{ v-\text{complex}}\\{l_v\equiv 0\,(r);\;l_{\bar
v} = 0}\\{l_v + re_v > 0}\\{\alpha_{v} := 1 + k_v - l_v - re_v \ge
0}}}\RN_{L/F}(\p_v) ^{\alpha_{v} - 1}\;\;\;\;
\cdot\prod_{\substack{{ v-\text{real}}\\{e_v
> 0}\\{\beta_{v} := 1 + k_v - re_v \ge 0}}} \q_{v} ^{\beta_{v} -
1}.
\end{eqnarray*}
 Recall that $\n ^{*}$ denotes the part of $\n$ coprime to
$\m_{1}\h_{1}.$ It follows that 
\begin{eqnarray*}
\n ^{*}& =& \n_{1}\cdot\,
\left(\frac{\overline{\m}}{\overline{\m}_{0}\overline{\m}_{2}}\right)\cdot\,
\RN_{L/F}(\m_{2})\cdot \h_{2}\\
&&\cdot \prod_{\substack{{
v-\text{complex}}\\{l_v\equiv 0\,(r);\;l_{\bar v} = 0}\\{l_v +
re_v
> 0}\\{\alpha_{v} := 1 + k_v - l_v - re_v \ge 0}}} \RN_{L/F}(\p_v)
^{\alpha_{v} - 1}\;\;\; \cdot\prod_{\substack{{
v-\text{real}}\\{e_v
> 0}\\{\beta_{v} := 1 + k_v - re_v \ge 0}}} \q_{v} ^{\beta_{v} -
1}.
\end{eqnarray*}

Combining all these with $(4.26),$ we obtain \begin{footnotesize}
$$\hskip-6pt\begin{aligned} &\sum_{\substack{{\frak{m}\in
\CI_{L}(S)}\\{\m-\text{imaginary}}}}
  \frac{(\psi
\tilde{\rho})(\m)}{\RN_{L/\BQ}(\m) ^s}\; \sum_{\substack{{\h\in
\CI_{L}(S)}\\{\n \in \CI_F(S)}\\{\h-\text{real}}}} \frac{(\psi
\tilde{\rho})(\h)\,\overline{\rho(\frak{n})}\,\,\overline{\chi_{\m_{1}\h 
_{1}}
^{*}(\n^*)}\,G(\chi_{\m_{1} \h_{1}} ^{*})\,G_{0}(\n, \m
\h)}{\RN_{L/\BQ}(\h) ^s\,\RN_{F/\BQ}(\n) ^w}\\
&= \sum_{\substack{{\frak{m}\in
\CI_{L}(S)}\\{\m-\text{imaginary}}}}
\frac{\psi(\m)\tilde{\rho}(\m_{0})\,\overline{\chi_{\m_{1}} ^{*}
\left(\frac{\overline{\m}}{\overline{\m}_{0}}\right)}\,\RN_{L/ 
\BQ}(\m_{0})
^{w - \frac{1}{2}}}{\RN_{L/\BQ}(\m) ^{s + w - \frac{1}{2}}} \\
&\cdot\sum_{\h\in \CI_F(S)} \frac{(\psi
\rho)(\h)\,\rho(\h_{0})\,\RN_{F/\BQ}(\h_{0}) ^{w -
1}\,\chi_{\h_{1}} ^{*}(\m)\,\chi_{\h_{1}} ^{*}(\m_{0})
^{-1}G(\chi_{\m_{1}\h_{1}} ^{*})}{\RN_{F/\BQ}(\h) ^{2s + w -
1}}\prod_{\substack{{ v}\\{\ord_{v}(\RN_{L/F}(\m_{1})) >
0}\\{\ord_{v}(\h_{2})
> 0}}} (1 - q_v ^{-1})\\&\cdot \prod_{\substack{{  
  v}\\{\ord_{v}(\RN_{L/F}(\m_{2}))
> 0}\\{\ord_{v}(\h)=0
}}} \Big [\,-\;(\chi_{\m_{1}} ^{*}\, \rho)(\pi_{v})\,q_{v} ^{w -
1}\, +\; (1 - q_v ^{-1})\;\cdot \sum_{\alpha_{v} \ge
0}(\,(\chi_{\m_{1}} ^{*}\,\rho) ^{-1}(\pi_{v})\,q_{v}
^{-w}\,)^{\alpha_{v}}\Big ]\\&\cdot \prod_{\substack{{
v}\\{\ord_{v}(\RN_{L/F}(\m_{2})) > 0}\\{\ord_{v}(\h_{2}) > 0 }}}
\Big [\,-\;(\chi_{\m_{1}} ^{*}\, \rho)(\pi_{v})\,q_{v} ^{w - 1}(1
- q_v ^{-1})\; +\; (1 - q_v ^{-1}) ^{2}\,\cdot \sum_{\alpha_{v}
\ge 0}(\,(\chi_{\m_{1}} ^{*}\,\rho) ^{-1}(\pi_{v})\,q_{v}
^{-w}\,)^{\alpha_{v}}\Big ]\\&\cdot\prod_{\substack{{v-\text{split
in $L$}}\\{\ord_{v}(\RN_{L/F}(\m)) = 0}\\{\ord_{v}(\h_{2}) > 0}}}
\Big [\,(\chi_{\m_{1}} ^{*}\, \rho)(\pi_{v})\,q_{v} ^{w - 2}\, +\;
(1 - q_v ^{-1}) ^{2}\,\cdot \sum_{\alpha_{v} \ge
0}(\,(\chi_{\m_{1}} ^{*}\,\rho) ^{-1}(\pi_{v})\,q_{v}
^{-w}\,)^{\alpha_{v}}\Big ]  \end{aligned}$$
$$\begin{aligned}&\cdot
\prod_{\substack{{v-\text{inert in $L$}}\\{\ord_{v}(\h_{2}) > 0}}}
\Big [\,-\;(\chi_{\m_{1}} ^{*}\, \rho)(\pi_{v})\,q_{v} ^{w - 2}\,
+\; (1 - q_v ^{-2})\,\cdot \sum_{\beta_{v} \ge 0}(\,(\chi_{\m_{1}}
^{*}\,\rho) ^{-1}(\pi_{v})\,q_{v} ^{-w}\,)^{\beta_{v}}\Big ]
\\&\cdot \sum_{\substack{{\n_{1} \in
\CI_F(S)}\\{(\n_{1},\m\h)\,=\,1}}}\; 
\frac{\overline{\rho(\n_{1})}\,\,\overline{\chi_{\m_{1}}
^{*}(\n_{1})}}{\RN_{F/\BQ}(\n_{1}) ^w}.
\end{aligned}$$
\end{footnotesize}

\noindent Note that the last sum represents an incomplete Hecke
$L$-function. After evaluating the geometric series inside the
last four products, the missing Euler factors corresponding to
places of $F$ dividing $\RN_{L/F}(\m_2) \h_2$ can be incorporated.
Also, multiply and divide by the Euler factors corresponding to
places of $F$ dividing $\h_0,$ forcing in this way $L_{S}(w,\,
(\chi_{\m_{1}} ^{*}\,\rho) ^{-1})$ to appear.

Let $R_{\frak{c}} ^{+}$ be the subgroup of $R_{\frak{c}}$
generated by the images (in $R_{\frak{c}}$) of all real fractional
ideals of $L$ coprime to $S'.$ Let $\e'$ be a fixed element of
$R_{\frak{c}}$ which is the image of an imaginary ideal $\m\in
\CI_L(S).$ Replacing $\psi$ by $\psi\tau_{1}\tau_{2}$ with
$\tau_{1}$ and $\tau_{2}$ characters of $R_{\frak{c}}$ and
$R_{\frak{c}}/R_{\frak{c}} ^{+},$ respectively, and making a
standard linear combination, one can restrict the first two sums
over ideals $\m$ and $\h,$ for which the image of $\m_1$ in
$R_{\frak{c}}$ is $\e'$ modulo $R_{\frak{c}} ^{+}$ and the image
of $\m_1 \h_1$ is a fixed element $\e$ of $R_\frak{c}$.

Now, invoke the functional equation of $L(w,\, (\chi_{\m_{1}}
^{*}\,\rho) ^{-1}).$ It is well-known, see \cite{W1}, that the
incomplete Hecke $L$-function (defined over $F$)
$$L_{S}\left(w,\, (\chi_{\m_{1}} ^{*}\,\rho) ^{-1}\right)\,=\,
\prod_{v\not\in S}L_{v}\left(w,\, (\chi_{\m_{1}} ^{*}\,\rho)_{v}
^{-1}\right)\,=\,\prod_{v\not\in S}\,\big [1\, -\, (\chi_{\m_{1}}
^{*}\,\rho)_{v} ^{-1}(\pi_{v})\,q_{v} ^{-w} \big ] ^{-1}$$
satisfies the functional equation 
\begin{eqnarray*}
L_{S}\left(w,\, (\chi_{\m_{1}} ^{*}\,\rho)
^{-1}\right) &=&\varepsilon\left(w,\, (\chi_{\m_{1}} ^{*}\,\rho)
^{-1}\right)\cdot L_{S}\left(1 - w,\, \chi_{\m_{1}}
^{*}\,\rho\right)\\&&\cdot\prod_{v\in S_{\infty}}\frac{L_{v}\left(1
- w,\, \rho_{v} \right)}{L_{v}\left(w,\, \rho_{v}
^{-1}\right)}\cdot\prod_{v\in S'}\frac{L_{v}\left(1 - w,\,
(\chi_{\m_{1}} ^{*}\,\rho)_{v} \right)}{L_{v}\left(w,\,
(\chi_{\m_{1}} ^{*}\,\rho)_{v} ^{-1}\right)}.\end{eqnarray*}
Replace $\psi$ by $\psi\,\eta ^{-1},$ and combine the above
functional equation with Lemma $3.6.$ Here $\Re (s)$ is taken
sufficiently large to ensure convergence. Using the
Fisher-Friedberg extension of the reciprocity law \cite{FF1}, one
can see that
$$\overline{\chi_{\m_{1}} ^{*} (\overline{\m})}\,\chi_{\h_{1}}
^{*}(\m)\,=\,C_{\e,\, \widehat{\e'}}'\cdot \chi_{\m}
^{*}(\h_{1}),$$ where $C_{\e,\, \widehat{\e'}}'$ is a constant
depending on just $\e$ and the class $\widehat{\e'}$ in
$R_{\frak{c}}/R_{\frak{c}} ^{+}.$ Also, note that $$\prod_{v\in
S'} \Bigl(1 - \rho ^{-r}(\pi_v)\,q_v ^{-rw}\Bigr)
^{-1}\cdot\;\,\frac{L_{v}\left(1 - w,\, (\chi_{\m_{1}}
^{*}\,\rho)_{v} \right)}{L_{v}\left(w,\, (\chi_{\m_{1}}
^{*}\,\rho)_{v} ^{-1}\right)}$$ is the inverse of a polynomial in
the variables $q_v ^{w},\, q_v ^{-w}$ corresponding to places
$v\in S'$ of the totally real field $F.$ The characters involved
in its coefficients are trivial on real ideals. Now, the
functional equation $(3.8)$ immediately follows, after we  replace
$\psi$ with $\psi \tau,$ where $\tau$ ranges over a finite set of
id\'ele class characters unramified outside $S$ and orders
dividing $r,$ and make a combination such that the above product
over $v\in S'$ disappears.

Starting from the definition of
$$\prod_{v\in S'} \Bigl(1 - \rho ^{r}(\pi_v)\,q_v ^{rw - r}\Bigr)
^{-1}\cdot\;\;\widetilde{Z}(s + w - {\scriptstyle \frac{1}{2}}, 1
- w; \psi; \rho),$$ one can easily check $(3.9)$ by reversing the
above argument.
\hfill\qed

\end{document}